\documentclass{elsarticle}
\usepackage[utf8]{inputenc}
\usepackage{lineno,hyperref}
\modulolinenumbers[5]
\usepackage{booktabs,makecell}
\usepackage{multicol}
\usepackage{multirow}
\usepackage{amsmath,stmaryrd}
\usepackage{graphicx}
\usepackage{subfigure}
\usepackage{bm}
\usepackage{enumerate}
\usepackage{pgfplots}
\usepackage{tikz}
\usepackage{algorithm}
\usepackage[noend]{algpseudocode}
\biboptions{numbers,sort&compress}
\usepackage{color,xcolor,colortbl}
\usepackage{footnote}
\usepackage{url}
\usepackage{array}
\usepackage{etex}
\usepackage{amsfonts}

\begin{document}

\begin{frontmatter}

		\title{The Cross Products of $M$ Vectors in $N$-dimensional Spaces and Their Geometric Significance}
		
        \author[rvt]{Chengshen Xu\corref{cor1}}
        \cortext[cor1]{Corresponding author: xcssgzs@126.com}
        \address[rvt]{Autohome Inc., 10th Floor Tower B, No. 3 Dan Ling Street Haidian District, Beijing, China}

\begin{abstract}

In textbooks and historical literature, the cross product has been defined only in $2$-dimensional and $3$-dimensional Euclidean spaces and the cross product of only two vectors has been defined only in the high dimensional Euclidean space whose metric matrix is the unit matrix. Nobody has given a universal definition for any number of vectors in high dimensional spaces whose metric matrices are the unit matrices. In fact, we can also define the cross product of $m$ vectors in an $n$-dimensional space, where $n$ and $m$ can take any positive integers larger than 1 and $m$ must not be larger than $n$. In this paper, we give the definition of the cross product of $m$ vectors in $n$-dimensional spaces whose metric matrices are any real symmetric or Hermitian matrices, and put forward two theorems related to matrices, so as to perfectly explain the geometric meaning of the cross product of vectors in high dimensional spaces. Furthermore, the length of the cross product represents the $m$-dimensional volume of the parallel polyhedron spanned by the $m$ vectors, and the absolute value of each component of the cross product represents each component of the volume in different directions. Specially, in the high dimensional Euclidean and unitary space where the metric matrix is the unit matrix, this decomposition of the volume still satisfies the Pythagorean theorem, and the cross product of $n$ vectors in an $n$-dimensional space is the determinant of the square matrix which is formed with these $n$ vectors as row or column vectors. We also explain the geometric meaning of the determinants of the metric matrices and their sub-matrices, which is also useful for understanding the invariant volume elements in high dimensional spaces and their subspaces in the differential geometry.
			
\end{abstract}
		
		
\end{frontmatter}
	
\section{Introduction}

In this paper, we use ``$\times$" to represent the cross product and use ``$\cdot$" to represent the inner product of vectors.

In textbooks and historical literature, the cross product is defined only in $2$-dimensional and $3$-dimensional Euclidean spaces and the cross product of only two vectors is defined in high dimensional spaces which are only Euclidean \cite{plebanski1988notes, deschapelles1993generalized, massey1983cross}. In these paper it is concluded that the product of $m$ vectors cannot be defined in high dimensional Euclidean spaces where $n>3$ and $2<m<n-1$ because the map
\begin{eqnarray}
\underbrace{\mathbb{R}^{n}\times\mathbb{R}^{n}\times\cdots\times\mathbb{R}^{n}}_{the ~total ~number ~is ~ m}\rightarrow\mathbb{R}^{n}
\label{1}
\end{eqnarray}
cannot satisfy the following conditions
\begin{eqnarray}
(X_{1}\times X_{2}\times\cdots\times X_{m})\cdot X_{i}=0
\end{eqnarray}
when $1\leq i\leq m$ and $2<m<n-1$ in high dimensional Euclidean spaces where $n>3$. Here $\mathbb{R}$ represents the real number field and $X_{i}$ represents an $n$-dimensional vector, i.e., an $n$-element array in field $\mathbb{R}$.

But we will know we should not define the cross product as Eq.(\ref{1}) when we can understand the geometric meaning of cross product clearly. In this paper, we define the cross product of $m$ vector in $n$-dimensional spaces over field $\mathbb{P}$ where $\mathbb{P}$ can be any number fields such as the real number field and the complex number field. What's more, the metric matrices of the spaces can be any real symmetric or Hermitian matrices. Unlike paper \cite{plebanski1988notes, deschapelles1993generalized, massey1983cross}, we define the cross product as follows:
\begin{eqnarray}
\underbrace{\mathbb{P}^{n}\times\mathbb{P}^{n}\times\cdots\times\mathbb{P}^{n}}_{the ~total ~number ~is ~ m}\rightarrow\mathbb{P}^{C^{m}_{n}}
\label{2}
\end{eqnarray}
where $m\leq n$ and $C^{m}_{n}=\frac{n!}{m!(n-m)!}$. The specific definition will be explained later. The absolute value of each component of the cross product represents each component of the volume of the parallel polyhedron spanned by the $m$ vectors in different directions, which will be proved later in this paper. When the metric matrix $G$ is the unit matrix
\begin{eqnarray}
G=I_{n}=
\left(
  \begin{array}{cccc}
    1 & 0 & \cdots & 0 \\
    0 & 1 & \cdots & 0 \\
    \vdots & \vdots &  & \vdots \\
    0 & 0 & \cdots & 1 \\
  \end{array}
\right)
\end{eqnarray}
the decomposition of the volume still satisfies the Pythagorean theorem, which will also be proved.

In Section 2, we define the cross product of $2$ vectors in the $2$-dimensional Euclidean spaces firstly and
prove that the absolute value of the cross product represents the volume of the parallelogram spanned by the $2$ vectors. Then we generalize the definition to the case $n$ vectors in the $n$-dimensional Euclidean spaces where the cross product is the determinant of the square matrix which is formed with these n vectors as row or column vectors and prove the absolute value of the determinant represents the volume of the parallel polyhedron spanned by the $n$ vectors. In that section, we only discuss for the case in which the metric matrix is the unit matrix and $\mathbb{P}$ is the real number field $\mathbb{R}$.

In Section 3, we put forward two theorems about matrices, based on which we can prove the rationality of the cross product definition and perfectly explain its geometric meaning.

In Section 4, we generalize the definition to the case $m$ vectors in $n$-dimensional Euclidean and unitary spaces based on the two theorems. The result vector has $C_{n}^{m}$ components. We explain the geometric meaning of its each component as mentioned above.


In Section 5, we generalize all the conclusions naturally to high dimensional spaces where the metric matrices can be any real symmetric or Hermitian matrices and $\mathbb{P}$ is the complex number field $\mathbb{C}$. These conclusions are also useful for understanding invariant volume elements in high dimensional spaces and their subspaces in the differential geometry.

Therefore, we suppose $G=I_{n}$ and $\mathbb{P}$ is the real number field $\mathbb{R}$ in Section 2 and $G=I_{n}$ and $\mathbb{P}$ is the complex number field $\mathbb{C}$ in Section 4.

\section{The Cross Product of $N$ Vectors in the $N$-dimensional Space Where the Metric Matrix Is the Unit Matrix}

In the $2$-dimensional Euclidean space with $G=I_{2}$, the area of the parallelogram spanned by $X_{1}=(x_{11}, x_{21})^{T}\in\mathbb{R}^{2}$ and $X_{2}=(x_{12}, x_{22})^{T}\in\mathbb{R}^{2}$, where the subscript ``$T$" represents the transpose operation, is
\begin{eqnarray}
S=||X_{1}||X_{2}|\sin\theta|=|\sqrt{x_{11}^{2}+x_{21}^{2}}\sqrt{x_{11}^{2}+x_{21}^{2}}\sin\theta|
\end{eqnarray}
where $|\cdot|$ represents the length of a vector or the absolute value of a number, and $\theta=<X_{1}, X_{2}>$ is the angle in the counterclockwise direction from $X_{1}$ to $X_{2}$. Now we define the cross product of $2$ vectors in the $2$-dimensional Euclidean space as
\begin{eqnarray}
X_{1}\times X_{2}=|X_{1}||X_{2}|\sin<X_{1}, X_{2}>.
\end{eqnarray}
Obviously, this definition satisfies the distributive law of multiplication when we make parallelogram factorization of vectors, i.e., $\forall ~Y_{1}, Y_{2}, Y_{3}, Y_{4}\in\mathbb{R}^{2}$, we have
\begin{eqnarray}
(Y_{1}+Y_{2})\times(Y_{3}\times Y_{4})=Y_{1}\times Y_{3}+Y_{1}\times Y_{4}+Y_{2}\times Y_{3}+Y_{2}\times Y_{4}.
\end{eqnarray}
Therefore we have
\begin{eqnarray}
&&
X_{1}\times X_{2}=[(x_{11}, 0)^{T}+(0, x_{21})^{T}]\times[(x_{12}, 0)^{T}+(0, x_{22})^{T}]
\nonumber\\
&&
(x_{11}, 0)^{T}\times(x_{12}, 0)^{T}+(x_{11}, 0)^{T}\times(0, x_{22})^{T}+(0, x_{21})^{T}\times(x_{12}, 0)^{T}
\nonumber\\
&&
+(0, x_{21})^{T}\times (0, x_{22})^{T}
\nonumber\\
&&
=x_{11}x_{12}\sin0+x_{11}x_{22}\sin\frac{\pi}{2}+x_{12}x_{21}\sin(-\frac{\pi}{2})+x_{21}x_{22}\sin0
\nonumber\\
&&
=x_{11}x_{22}-x_{21}x_{12}
=
\left|
  \begin{array}{cc}
    x_{11} & x_{12} \\
    x_{21} & x_{22} \\
  \end{array}
\right|,
\end{eqnarray}
So the absolute value of the determinant is the area of the parallelogram for $n=2$.

Obviously, the cross product also satisfies antisymmetry, i.e., $\forall ~Y_{1}, Y_{2}\in\mathbb{R}^{2}$,
\begin{eqnarray}
Y_{1}\times Y_{2}=-Y_{2}\times Y_{1}
\end{eqnarray}

Similarly, using the distributive law and antisymmetry, we can find the volume of a parallel polyhedron spanned by
\begin{eqnarray}
\begin{array}{c}
  X_{1}=(x_{11}, x_{21}, \cdots, x_{n1})^{T}\in\mathbb{R}^{n} \\
  X_{2}=(x_{12}, x_{22}, \cdots, x_{n2})^{T}\in\mathbb{R}^{n} \\
  \vdots \\
  X_{n}=(x_{1n}, x_{2n}, \cdots, x_{nn})^{T}\in\mathbb{R}^{n}
\end{array}
\label{vectorsR}
\end{eqnarray}
in the $n$-dimensional Euclidean space where $g=I_{n}$ is
\begin{eqnarray}
V=|\sum\limits_{i_{1}, i_{2}, \cdots, i_{n}}(-1)^{\tau(i_{1}, i_{2}, \cdots, i_{n})}x_{1i_{1}}x_{2i_{2}}\cdots x_{ni_{n}}|
=
\left|\left|
  \begin{array}{cccc}
    x_{11} & x_{12} & \cdots & x_{1n} \\
    x_{21} & x_{22} & \cdots & x_{2n} \\
    \vdots & \vdots &  & \vdots \\
    x_{n1} & x_{n2} & \cdots & x_{nn} \\
  \end{array}
\right|\right|
\end{eqnarray}
where $i_{1}, i_{2}, \cdots, i_{n}$ is a permutation of $1, 2, \cdots, n$ and $\tau(i_{1}, i_{2}, \cdots, i_{n})$ is the inversion number of the permutation. So we define the cross product of $n$ vectors in the $n$-dimensional Euclidean space as
\begin{eqnarray}
X_{1}\times X_{2}\times\cdots\times X_{n}
=
\left|
  \begin{array}{cccc}
    x_{11} & x_{12} & \cdots & x_{1n} \\
    x_{21} & x_{22} & \cdots & x_{2n} \\
    \vdots & \vdots &  & \vdots \\
    x_{n1} & x_{n2} & \cdots & x_{nn} \\
  \end{array}
\right|
\end{eqnarray}
whose absolute value is the volume of the parallel polyhedron spanned by the vectors (\ref{vectorsR}) in the $n$-dimensional Euclidean space where $G=I_{n}$.

\section{Two Important Theorems Found by Us}

In this section we prove two theorems about matrices in order to generalize the cross product in next two section.

Theorem 1. $\forall ~$ matrix A with $m$ rows and $n$ columns
\begin{eqnarray}
A=
\left(
  \begin{array}{cccc}
    a_{11} & a_{12} & \cdots & a_{1n} \\
    a_{21} & a_{22} & \cdots & a_{2n} \\
    \vdots & \vdots &  & \vdots \\
    a_{m1} & a_{m2} & \cdots & a_{mn} \\
  \end{array}
\right)
\end{eqnarray}
and $\forall ~$ matrix B with $n$ rows and $m$ columns
\begin{eqnarray}
B=
\left(
  \begin{array}{cccc}
    b_{11} & b_{12} & \cdots & b_{1m} \\
    b_{21} & b_{22} & \cdots & b_{2m} \\
    \vdots & \vdots &  & \vdots \\
    b_{n1} & b_{n2} & \cdots & b_{nm} \\
  \end{array}
\right)
\end{eqnarray}
we have
\begin{eqnarray}
|AB|=
\left\{
  \begin{array}{ll}
    0, & m>n; \\
    \left(
      \begin{array}{cccc}
        |A_{1}| & |A_{2}| & \cdots & |A_{C_{n}^{m}}| \\
      \end{array}
    \right)
    \left(
      \begin{array}{c}
        |B_{1}| \\
        |B_{2}| \\
        \vdots \\
        |B_{C_{n}^{m}}| \\
      \end{array}
    \right)
, & m\leq n.
  \end{array}
\right.
\end{eqnarray}
where the columns of Matrix $A_{i}$ are $m$ different columns of $A$ with the column indexes from small to large and the rows of Matrix $B_{i}$ are $m$ different rows of $B$  with the row indexes from small to large. The selection for the columns of $A_{i}$ corresponds to that for the rows of $B_{i}$. Different values of $i$ correspond to different combinations, so the total number is $C_{n}^{m}$.

Proof. When $m>n$, the conclusion is clearly valid. When $m\leq n$, we have
\begin{eqnarray}
&&
\left(
  \begin{array}{cccc}
    a_{11} & a_{12} & \cdots & a_{1n} \\
    a_{21} & a_{22} & \cdots & a_{2n} \\
    \vdots & \vdots &        & \vdots \\
    a_{m1} & a_{m2} & \cdots & a_{mn} \\
  \end{array}
\right)
\left(
  \begin{array}{cccc}
    b_{11} & b_{12} & \cdots & b_{1m} \\
    b_{21} & b_{22} & \cdots & b_{2m} \\
    \vdots & \vdots &        & \vdots \\
    b_{n1} & b_{n2} & \cdots & b_{nm} \\
  \end{array}
\right)
\nonumber\\
&&
=
\left(
  \begin{array}{cccc}
    \sum\limits^{n}_{i=1}a_{1i}b_{i1} & \sum\limits^{n}_{i=1}a_{1i}b_{i2} & \cdots & \sum\limits^{n}_{i=1}a_{1i}b_{im} \\
    \sum\limits^{n}_{i=1}a_{2i}b_{i1} & \sum\limits^{n}_{i=1}a_{2i}b_{i2} & \cdots & \sum\limits^{n}_{i=1}a_{2i}b_{im} \\
    \vdots & \vdots &        & \vdots \\
    \sum\limits^{n}_{i=1}a_{mi}b_{i1} & \sum\limits^{n}_{i=1}a_{mi}b_{i2} & \cdots & \sum\limits^{n}_{i=1}a_{mi}b_{im} \\
  \end{array}
\right)
\nonumber\\
&&
=
\sum\limits_{j_{1},j_{2}, \cdots, j_{m}}(-1)^{\tau(j_{1}, j_{2}, \cdots, j_{m})}(\sum\limits_{i=1}^{n}a_{1i}b_{ij_{1}})(\sum\limits_{i=1}^{n}a_{2i}b_{ij_{2}})
\cdots(\sum\limits_{i=1}^{n}a_{mi}b_{ij_{m}})
\nonumber\\
&&
=
\sum\limits_{j_{1},j_{2}, \cdots, j_{m}}(-1)^{\tau(j_{1}, j_{2}, \cdots, j_{m})}\sum\limits_{l_{1}, l_{2}, \cdots, l_{m}}a_{1l_{1}}b_{l_{1}j_{1}}a_{2l_{2}}b_{l_{2}j_{2}}
\cdots a_{ml_{m}}b_{l_{m}j_{m}}
\nonumber\\
&&
=\sum\limits^{C_{n}^{m}}_{i=1}|A_{i}B_{i}|=\sum\limits^{C_{n}^{m}}_{i=1}|A_{i}||B_{i}|
\nonumber\\
&&
=
\left(
      \begin{array}{cccc}
        |A_{1}| & |A_{2}| & \cdots & |A_{C_{n}^{m}}| \\
      \end{array}
    \right)
    \left(
      \begin{array}{c}
        |B_{1}| \\
        |B_{2}| \\
        \vdots \\
        |B_{C_{n}^{m}}| \\
      \end{array}
    \right),
\end{eqnarray}
where $j_{1}, j_{2}, \cdots, j_{m}$ is a permutation of $1, 2, \cdots, m$, $\tau(j_{1}, j_{2}, \cdots, j_{m})$ is the inversion number of the permutation, $l_{1}, l_{2}, \cdots, l_{m}$ is a permutation of $m$ elements out of $1, 2, \cdots n$. At the third step we have used the antisymmetry to cancel out a lot of terms.

Theorem 2. $\forall ~$ matrix A with $m$ rows and $n$ columns
\begin{eqnarray}
A=
\left(
  \begin{array}{cccc}
    a_{11} & a_{12} & \cdots & a_{1n} \\
    a_{21} & a_{22} & \cdots & a_{2n} \\
    \vdots & \vdots &  & \vdots \\
    a_{m1} & a_{m2} & \cdots & a_{mn} \\
  \end{array}
\right)
\end{eqnarray}
and $\forall ~$ matrix B with $n$ rows and $s$ columns
\begin{eqnarray}
B=
\left(
  \begin{array}{cccc}
    b_{11} & b_{12} & \cdots & b_{1s} \\
    b_{21} & b_{22} & \cdots & b_{2s} \\
    \vdots & \vdots &  & \vdots \\
    b_{n1} & b_{n2} & \cdots & b_{ns} \\
  \end{array}
\right)
\end{eqnarray}
and
\begin{eqnarray}
C=AB
\end{eqnarray}
with $m$ rows and $s$ columns, we have
\begin{eqnarray}
&&
\left(
  \begin{array}{cccc}
    |C_{11}|         & |C_{12}|         & \cdots & |C_{1C^{k}_{s}}| \\
    |C_{21}|         & |C_{22}|         & \cdots & |C_{2C^{k}_{s}}| \\
    \vdots           & \vdots           &        & \vdots \\
    |C_{C^{k}_{m}1}| & |C_{C^{k}_{m}2}| & \cdots & |C_{C^{k}_{m}C^{k}_{s}}| \\
  \end{array}
\right)
=
\left(
  \begin{array}{cccc}
    |A_{11}|         & |A_{12}|         & \cdots  & |A_{1C^{k}_{n}}| \\
    |A_{21}|         & |A_{22}|         & \cdots  & |A_{2C^{k}_{n}}| \\
    \vdots           & \vdots           &         & \vdots \\
    |A_{C^{k}_{m}1}| & |A_{C^{k}_{m}2}| & \cdots  & |A_{C^{k}_{m}C^{k}_{n}}| \\
  \end{array}
\right)
\nonumber\\
&&
\left(
  \begin{array}{cccc}
    |B_{11}|         & |B_{12}|         & \cdots & |B_{1C^{k}_{s}}| \\
    |B_{21}|         & |B_{22}|         & \cdots & |B_{2C^{k}_{s}}| \\
    \vdots           & \vdots           &        & \vdots \\
    |B_{C^{k}_{n}1}| & |B_{C^{k}_{n}2}| & \cdots & |B_{C^{k}_{n}C^{k}_{s}}| \\
  \end{array}
\right),
\end{eqnarray}
where $C_{ij}$ is the sub-matrix generated by the elements from $k$ different rows and $k$ different columns of $C$ with the row and column indexes from small to large, $A_{ij}$ is the sub-matrix generated by the elements from $k$ different rows and $k$ different columns of $A$ with the row and column indexes from small to large and $B_{ij}$ is the sub-matrix generated by the elements from $k$ different rows and $k$ different columns of $B$ with the row and column indexes from small to large. The selection for the row indexes of $C_{ij}$ corresponds to that of $A_{il}$, the selection for the column indexes of $C_{ij}$ corresponds to that of $B_{lj}$ and the selection for the column indexes of $A_{il}$ corresponds to the selection for the row indexes of $A_{lj}$. Different values of $i$, $j$ or $l$ correspond to different combinations, so the total numbers are $C_{m}^{k}$, $C_{n}^{k}$ and $C_{s}^{k}$ respectively.

Proof. $\forall ~1\leq i\leq C_{m}^{k}, 1\leq j\leq C_{s}^{k}$, According to theorem $1$, we have
\begin{eqnarray}
|C_{ij}|=\sum\limits_{l=1}^{C_{n}^{k}}|A_{il}||B_{lj}|
\end{eqnarray}
so the conclusion is valid.

\section{The Cross Product of $M$ Vectors in the $N$-dimensional Euclidean and Unitary Space Where the Metric Matrix Is the Unit Matrix}

According to the two theorems in Section 3, we can define the cross product of $m$ vectors in the $n$-dimensional Euclidean space where the metric matrix is the unit matrix. Supposing
\begin{eqnarray}
\begin{array}{c}
  X_{1}=(x_{11}, x_{21}, \cdots, x_{n1})^{T}\in\mathbb{R}^{n} \\
  X_{2}=(x_{12}, x_{22}, \cdots, x_{n2})^{T}\in\mathbb{R}^{n} \\
  \vdots \\
  X_{m}=(x_{1m}, x_{2m}, \cdots, x_{nm})^{T}\in\mathbb{R}^{n}
\end{array}
\end{eqnarray}
and $2\leq m\leq n$, we make
\begin{eqnarray}
X=
\left(
  \begin{array}{cccc}
    X_{1} & X_{2} & \cdots & X_{m} \\
  \end{array}
\right)
=
\left(
  \begin{array}{cccc}
    x_{11} & x_{12} & \cdots & x_{1m} \\
    x_{21} & x_{22} & \cdots & x_{2m} \\
    \vdots & \vdots &  & \vdots\\
    x_{n1} & x_{n2} & \cdots & x_{nm}
  \end{array}
\right)
\label{equation1}
\end{eqnarray}
and we define
\begin{eqnarray}
\tilde{X}
=
\left(
  \begin{array}{c}
    |X^{1}| \\
    |X^{2}| \\
    \vdots \\
    |X^{C_{n}^{m}}|
  \end{array}
\right),
\label{equation2}
\end{eqnarray}
where the rows of Matrix $X^{i}$ are $m$ different rows of $X$ with the row indexes from small to large. Different values of $i$ correspond to different combinations, so the total number is $C_{n}^{m}$. And $\forall$ Matric $A$ with $n$ rows and $n$ columns, we define
\begin{eqnarray}
\tilde{A}(m)
=
\left(
  \begin{array}{cccc}
    |A_{11}|         & |A_{12}|         & \cdots & |A_{1C^{m}_{n}}| \\
    |A_{21}|         & |A_{22}|         & \cdots & |A_{2C^{m}_{n}}| \\
    \vdots           & \vdots           &        & \vdots \\
    |A_{C^{m}_{n}1}| & |A_{C^{m}_{n}2}| & \cdots & |A_{C^{m}_{n}C^{m}_{n}}| \\
  \end{array}
\right),
\end{eqnarray}
where $A_{ij}$ is the sub-matrix generated by the elements from $m$ different rows and $m$ different columns of $A$ with the row and column indexes from small to large and the row indexes and their sorting order of $A_{li}$ are consistent with the column indexes and their sorting order of $A_{jl}$. Different values $i$, $j$ or $l$ correspond to different combinations, so the total number are $C_{n}^{m}$.

Now we define the cross product as
\begin{eqnarray}
X_{1}\times X_{2}\times\cdots\times X_{m}=\tilde{X}.
\label{equation3}
\end{eqnarray}
Obviously, it satisfies the distributive law of multiplication and the antisymmetry, and this cross product is zero vector when $X_{1}, X_{2}, \cdots, X_{m}$ linearly dependent.

From Section 2, we know the absolute value of the determinant of a matrix represents the $n$-dimensional volume of the parallel polyhedron spanned by the $n$ row or column vectors of the matrix, so the absolute value of each component of the vector $\tilde{X}$ represents each component of the volume of the parallel polyhedron, which is spanned by $X_{1}, X_{2}, \cdots, X_{m}$, in each $m$-dimensional subspace generated by $m$ basis vectors. In the process of the orthogonal decomposition in the $C_{n}^{m}$ directions, each decomposition satisfies the Pythagorean theorem from the geometric aspect. So the volume of the parallel polyhedron is
\begin{eqnarray}
V=|\tilde{X}|=\sqrt{\sum\limits_{i=1}^{C_{n}^{m}}||X^{i}||^{2}},
\end{eqnarray}
which we will prove right now from the algebraic aspect. Here the inner $|\cdot|$ represents the determinant and the outer $|\cdot|$ represents the absolute value.

In fact, we can choose a new set of orthogonal normalized basis vectors in the subspace generated by $X_{1}, X_{2}, \cdots, X_{m}$ and its orthogonal complement space or equivalently we can rotate the space to make $X_{1}, X_{2}, \cdots, X_{n}$ into the first $m$-dimensional subspace corresponding to $X^{1}$, i.e., there is a rotation matrix $S$ (i.e. a orthogonal matrix whose determinant is equal to $1$) that satisfies
\begin{eqnarray}
Y
=SX
=
S
\left(
  \begin{array}{cccc}
    X_{1} & X_{2} & \cdots & X_{m} \\
  \end{array}
\right)
=
\left(
  \begin{array}{cccc}
    y_{11} & y_{12} & \cdots & y_{1m} \\
    y_{21} & y_{22} & \cdots & y_{2m} \\
    \vdots & \vdots &   & \vdots \\
    y_{m1} & y_{m2} & \cdots & y_{mm} \\
    0 & 0 & \cdots & 0 \\
    0 & 0 & \cdots & 0 \\
    \vdots & \vdots &   & \vdots \\
    0 & 0 & \cdots & 0 \\
  \end{array}
\right).
\label{equation4}
\end{eqnarray}
Therefor according to Theorem 2 we have
\begin{eqnarray}
\tilde{Y}
=
\left(
  \begin{array}{c}
    |Y^{1}| \\
    0 \\
    0 \\
    0 \\
  \end{array}
\right)
=\tilde{S}(m)\tilde{X}
\end{eqnarray}
where
\begin{eqnarray}
Y^{1}
=
\left(
  \begin{array}{cccc}
    y_{11} & y_{12} & \cdots & y_{1m} \\
    y_{21} & y_{22} & \cdots & y_{2m} \\
    \vdots & \vdots &   & \vdots \\
    y_{m1} & y_{m2} & \cdots & y_{mm} \\
  \end{array}
\right),
\end{eqnarray}
so the volume of the parallel polyhedron is
\begin{eqnarray}
&&
V=||Y_{1}||=|\tilde{Y}|=\sqrt{\tilde{Y}\cdot\tilde{Y}}=\sqrt{\tilde{Y}^{T}\tilde{Y}}
=\sqrt{(\tilde{S}(m)\tilde{X})^{T}\tilde{S}(m)\tilde{X}}
\nonumber\\
&&
=\sqrt{\tilde{X}^{T}\tilde{S}(m)^{T}\tilde{S}(m)\tilde{X}}
=\sqrt{\tilde{X}^{T}\tilde{I_{n}}(m)\tilde{X}}
=\sqrt{\tilde{X}^{T}I_{C_{n}^{m}}\tilde{X}}
\nonumber\\
&&
=\sqrt{\tilde{X}^{T}\tilde{X}}=\sqrt{|X^{T}X|}=|\tilde{X}|.
\end{eqnarray}
It means that the volume of the parallel polyhedron, which is spanned by $X_{1}, X_{2}, \cdots, X_{m}$, is equal to the determinant of $X^{T}X$ to the power of one second as well as the length of $\tilde{X}$

So far we have proved that the length of the cross product of $m$ vectors in the $n$-dimensional Euclidean space whose metric matrix is the unit matrix we define represents the volume of the parallel polyhedron spanned by the $m$ vectors.

In the $n$-dimensional unitary space whose metric matrix is $I_{n}$, we have the similar conclusion. Supposing
\begin{eqnarray}
\begin{array}{c}
  X_{1}=(x_{11}, x_{21}, \cdots, x_{n1})^{T}\in\mathbb{C}^{n} \\
  X_{2}=(x_{12}, x_{22}, \cdots, x_{n2})^{T}\in\mathbb{C}^{n} \\
  \vdots \\
  X_{m}=(x_{1m}, x_{2m}, \cdots, x_{nm})^{T}\in\mathbb{C}^{n}
\end{array},
\label{equation5}
\end{eqnarray}
$2\leq m\leq n$, Eq.(\ref{equation1}) and Eq.(\ref{equation2}), we can also define the cross product as Eq.(\ref{equation3}). Similarly, we can choose a new set of orthogonal normalized basis vectors in the subspace generated by $X_{1}, X_{2}, \cdots, X_{m}$ and its orthogonal complement space or equivalently we can rotate the space to make $X_{1}, X_{2}, \cdots, X_{n}$ into the first $m$-dimensional subspace corresponding to $X^{1}$, i.e., there is a unitary matrix $S$ whose determinant is equal to $1$ that satisfies Eq.(\ref{equation4}). So we have
\begin{eqnarray}
&&
||Y_{1}||=|\tilde{Y}|=\sqrt{\tilde{Y}\cdot\tilde{Y}}=\sqrt{\tilde{Y}^{H}\tilde{Y}}
=\sqrt{(\tilde{S}\tilde{X})^{H}\tilde{S}\tilde{X}}
\nonumber\\
&&
=\sqrt{\tilde{X}^{H}\tilde{S}^{H}\tilde{S}\tilde{X}}
=\sqrt{\tilde{X}^{H}\tilde{I_{n}}(m)\tilde{X}}
=\sqrt{\tilde{X}^{H}I_{C_{n}^{m}}\tilde{X}}
\nonumber\\
&&
=\sqrt{\tilde{X}^{H}\tilde{X}}=\sqrt{|X^{H}X|}=|\tilde{X}|
\end{eqnarray}
where the superscript ``$H$" represents the transpose complex conjugate operation and $|\cdot|$ represents the absolute value for a complex number and the determinant for a matrix.

\section{The Cross Product of $M$ Vectors in $N$-dimensional Spaces Where the Metric Matrices Are Any Matrices}

A metric matrix in an $n$-dimensional space can be viewed as inner products of the pairs of the basis vectors $\bm{\varepsilon}_{1}, \bm{\varepsilon}_{2}, \cdots, \bm{\varepsilon}_{n}$
\begin{eqnarray}
G=
\left(
  \begin{array}{cccc}
    \bm{\varepsilon}_{1}\cdot\bm{\varepsilon}_{1} & \bm{\varepsilon}_{1}\cdot\bm{\varepsilon}_{2} & \cdots & \bm{\varepsilon}_{1}\cdot\bm{\varepsilon}_{n} \\
    \bm{\varepsilon}_{2}\cdot\bm{\varepsilon}_{1} & \bm{\varepsilon}_{2}\cdot\bm{\varepsilon}_{2} & \cdots & \bm{\varepsilon}_{2}\cdot\bm{\varepsilon}_{n} \\
    \vdots & \vdots &   & \vdots \\
    \bm{\varepsilon}_{n}\cdot\bm{\varepsilon}_{1} & \bm{\varepsilon}_{n}\cdot\bm{\varepsilon}_{2} & \cdots & \bm{\varepsilon}_{n}\cdot\bm{\varepsilon}_{n} \\
  \end{array}
\right)
\end{eqnarray}
which is a Hermitian matrix. If we view the basis vector as a column vector, we have
\begin{eqnarray}
&&
G=
\left(
  \begin{array}{c}
    \bm{\varepsilon}_{1}^{T} \\
    \bm{\varepsilon}_{2}^{T} \\
    \vdots \\
    \bm{\varepsilon}_{n}^{T} \\
  \end{array}
\right)
\left(
  \begin{array}{cccc}
    \bm{\varepsilon}_{1} & \bm{\varepsilon}_{2} & \cdots & \bm{\varepsilon}_{n} \\
  \end{array}
\right)
\nonumber\\
&&
\Longrightarrow
||G||=
\left|\left|
  \begin{array}{c}
    \bm{\varepsilon}_{1}^{T} \\
    \bm{\varepsilon}_{2}^{T} \\
    \vdots \\
    \bm{\varepsilon}_{n}^{T} \\
  \end{array}
\right|
\left|
  \begin{array}{cccc}
    \bm{\varepsilon}_{1} & \bm{\varepsilon}_{2} & \cdots & \bm{\varepsilon}_{n} \\
  \end{array}
\right|\right|
=
\left|
  \begin{array}{cccc}
    \bm{\varepsilon}_{1} & \bm{\varepsilon}_{2} & \cdots & \bm{\varepsilon}_{n} \\
  \end{array}
\right|^{2}.
\end{eqnarray}
According to Section 2, the determinant
$\left|
  \begin{array}{cccc}
    \bm{\varepsilon}_{1} & \bm{\varepsilon}_{2} & \cdots & \bm{\varepsilon}_{n} \\
  \end{array}
\right|$
represents the cross product of $\bm{\varepsilon}_{1}, \bm{\varepsilon}_{2}, \cdots, \bm{\varepsilon}_{n}$, so the absolute value of $|G|$ is the projection of the volume of the parallel polyhedron, which is spanned by $\bm{\varepsilon}_{1}, \bm{\varepsilon}_{2}, \cdots, \bm{\varepsilon}_{n}$, in the space generated by $\bm{\varepsilon}_{1}, \bm{\varepsilon}_{2}, \cdots, \bm{\varepsilon}_{n}$ multiplied by itself, i.e., the absolute value of $|G|$ is the volume of the parallel polyhedron spanned by $\bm{\varepsilon}_{1}, \bm{\varepsilon}_{2}, \cdots, \bm{\varepsilon}_{n}$ to the power of $2$. This is why the invariant volume element contains $\sqrt{||G||}$ in differential geometry. For
\begin{eqnarray}
\tilde{G}(m)
=
\left(
  \begin{array}{cccc}
    |G_{11}|         & |G_{12}|         & \cdots & |G_{1C^{m}_{n}}| \\
    |G_{21}|         & |G_{22}|         & \cdots & |G_{2C^{m}_{n}}| \\
    \vdots           & \vdots           &        & \vdots \\
    |G_{C^{m}_{n}1}| & |G_{C^{m}_{n}2}| & \cdots & |G_{C^{m}_{n}C^{m}_{n}}| \\
  \end{array}
\right)
\end{eqnarray}
we have
\begin{eqnarray}
&&
|G_{ij}|=
\left|
  \begin{array}{cccc}
    \bm{\varepsilon}_{i_{1}}\cdot\bm{\varepsilon}_{j_{1}} & \bm{\varepsilon}_{i_{1}}\cdot\bm{\varepsilon}_{j_{2}} & \cdots & \bm{\varepsilon}_{i_{1}}\cdot\bm{\varepsilon}_{i_{m}} \\
    \bm{\varepsilon}_{i_{2}}\cdot\bm{\varepsilon}_{j_{1}} & \bm{\varepsilon}_{i_{2}}\cdot\bm{\varepsilon}_{j_{2}} & \cdots & \bm{\varepsilon}_{i_{2}}\cdot\bm{\varepsilon}_{j_{m}} \\
    \vdots & \vdots &   & \vdots \\
    \bm{\varepsilon}_{i_{m}}\cdot\bm{\varepsilon}_{j_{1}} & \bm{\varepsilon}_{i_{m}}\cdot\bm{\varepsilon}_{j_{2}} & \cdots & \bm{\varepsilon}_{i_{m}}\cdot\bm{\varepsilon}_{j_{m}} \\
  \end{array}
\right|
\nonumber\\
&&
=
\left|\left|
  \begin{array}{c}
    \bm{\varepsilon}_{i_{1}}^{T} \\
    \bm{\varepsilon}_{i_{2}}^{T} \\
    \vdots \\
    \bm{\varepsilon}_{i_{m}}^{T} \\
  \end{array}
\right|
\left|
  \begin{array}{cccc}
    \bm{\varepsilon}_{j_{1}} & \bm{\varepsilon}_{j_{2}} & \cdots & \bm{\varepsilon}_{j_{m}} \\
  \end{array}
\right|\right|.
\end{eqnarray}
Therefore, the determinant $|G_{ij}|$ represents the projection of the cross product of $\bm{\varepsilon}_{i_{1}}, \bm{\varepsilon}_{i_{2}}, \cdots, \bm{\varepsilon}_{i_{m}}$ in the subspace generated by $\bm{\varepsilon}_{j_{1}}, \bm{\varepsilon}_{j_{2}}, \cdots, \bm{\varepsilon}_{j_{m}}$ multiplied by the cross product of  $\bm{\varepsilon}_{j_{1}}, \bm{\varepsilon}_{j_{2}}, \cdots, \bm{\varepsilon}_{j_{m}}$, i.e., the absolute value of $|G_{ij}|$ represents the projection of the volume of the parallel polyhedron spanned by $\bm{\varepsilon}_{i_{1}}, \bm{\varepsilon}_{i_{2}}, \cdots, \bm{\varepsilon}_{i_{m}}$ in the subspace generated by $\bm{\varepsilon}_{j_{1}}, \bm{\varepsilon}_{j_{2}}, \cdots, \bm{\varepsilon}_{j_{m}}$ multiplied by the volume of the parallel polyhedron spanned by $\bm{\varepsilon}_{j_{1}}, \bm{\varepsilon}_{j_{2}}, \cdots, \bm{\varepsilon}_{j_{m}}$, because of which we have
\begin{eqnarray}
\sqrt{\tilde{X}\cdot\tilde{X}}=\sqrt{\sum\limits_{i=1}^{C_{n}^{m}}\sum\limits_{j=1}^{C_{n}^{m}}
\overline{\tilde{X}}_{i}|G_{ij}|\tilde{X}_{j}}
=\sqrt{\tilde{X}^{H}\tilde{G}(m)\tilde{X}}=\sqrt{|X^{H}GX|}=|\tilde{X}|
\end{eqnarray}
where $\overline{\tilde{X}}_{i}$ represents the complex conjugate number of $\tilde{X}_{i}$. That means that our results can be generalized to the circumstances $m$ vectors in $n$-dimensional spaces whose metric matrices are any matrices if we define the cross product of Eq.(\ref{equation5}) as Eq.(\ref{equation3})

\section{Conclusions and Prospects}

In this paper, we have successfully generalized the definition of the cross product to any number of vectors in any dimensional spaces whose metric matrices are any Hermitian matrices. The conclusions we obtain are perfect whether they are at the algebraic aspect or the geometrical aspect. We put forward two theorems related to matrices, which makes us explain the geometric meaning of the cross product of vectors in high dimensional spaces. The length of the cross product represents the $m$-dimensional volume of the parallel polyhedron spanned by the $m$ vectors, and the absolute value of each component of the cross product represents each component of the volume in different directions. Specially, in the high dimensional Euclidean and unitary space when the metric matrix is the unit matrix, this decomposition of the volume still satisfies the Pythagorean theorem, and the cross product of $n$ vectors in an $n$-dimensional space is the determinant of the square matrix which is formed with these $n$ vectors as row or column vectors. We have also explained the geometric meaning of the determinants of the metric matrices and theirs sub-matrices, which is also useful for understanding the invariant volume elements in high dimensional spaces and their subspaces in differential geometry.

There are still issues that can be further discussed. For example, we can make more detailed discussions for the cases that the metric matrix is positive definite, positive semi-definite, negative definite, negative semi-definite, indefinite non-singular and indefinite singular respectively. For the case of the vector whose length is zero in an indefinite non-singular subspace, we cannot rotate the vector to a positive definite or a negative definite subspace such as light-like vectors in the special theory of relativity.

	
	
\section*{References}
    \bibliographystyle{elsart-num1}
	\bibliography{research}

\end{document}